\documentclass[reqno]{amsart}

\usepackage{amssymb}

\newcommand\Curr{\mathop {\fam 0 Cur}\nolimits}
\newcommand\Cend{\mathop {\fam 0 Cend}\nolimits}
\newcommand\Chom{\mathop {\fam 0 Chom}\nolimits}
\newcommand\gc{\mathop {\fam 0 gc}\nolimits}
\newcommand\End{\mathop {\fam 0 End}\nolimits}
\newcommand\Hom{\mathop {\fam 0 Hom}\nolimits}
\newcommand\Image{\mathop{\fam 0 Im}\nolimits}
\newcommand\Kernel{\mathop {\fam 0 Ker}\nolimits}
\newcommand\Rad{\mathop {\fam 0 Rad}\nolimits}

\newcommand\Idd{\mathop {\fam 0 Id}\nolimits}
\newcommand\Aw{{\mathrm A}}
\newcommand\rank{\mathop {\fam 0 rank}\nolimits}
\newcommand\diag{\mathop {\fam 0 diag}\nolimits}
\newcommand\oo[1]{\mathrel { \circ_{#1}}}

\newtheorem{theorem}{Theorem}
\newtheorem{lemma}[theorem]{Lemma}
\newtheorem{proposition}[theorem]{Proposition}
\newtheorem{corollary}[theorem]{Corollary}

\theoremstyle{definition}
\newtheorem{definition}[theorem]{Definition}
\newtheorem{example}[theorem]{Example}
\newtheorem{remark}[theorem]{Remark}

\numberwithin{equation}{section}
\numberwithin{theorem}{section}

\begin{document}

\title{On the Wedderburn
Principal Theorem in Conformal Algebras}

\thanks{Partially supported by RFBR 05--01--00230,
Complex Integration Program SB RAS (1.9)
and SB RAS grant for young researchers (Presidium SB RAS, act N.29 of January 26, 2006).
The author gratefully acknowledges the support of the Pierre Deligne fund
based on his 2004 Balzan prize in mathematics.}

\author{Pavel Kolesnikov}

\address{Sobolev Institute of Mathematics, Akad. Koptyug prosp., 4\\ Novosibirsk, 630090, Russia}

\email{pavelsk@math.nsc.ru}

\address{Korea Institute for Advanced Study, 207-43 Cheongnyangni 2-dong, Dongdaemun-gu \\
Seoul, 130722, Republic of Korea}

\begin{abstract}
We investigate an analogue of the Wedderburn principal theorem for
associative conformal algebras with finite faithful representations.
It is shown that the radical splitting property for an algebra of this kind
holds if the maximal semisimple factor of this algebra is unital,
but does not hold in general.
\end{abstract}



\maketitle

\section{Introduction}

One of the main points of the theory of finite-dimensional associative algebras
is the classical Wedderburn principal theorem. In a sketched
form, it could be stated as follows. If $A$ is an algebra over an
algebraically closed field $\Bbbk $, $N$ is a nilpotent ideal
of $A$, and if $A/N$ contains a semisimple finite-dimensional subalgebra
$\overline S$, then the preimage of $\overline{S}$ in $A$ contains a subalgebra
$S\simeq \overline{S}$. There are many generalizations of this theorem
for various classes of associative
(see, e.g., \cite{Az,Fe,W,Ze} and references therein)
and nonassociative (e.g., \cite{Al,E,Sc1,Sc2,SvO,Rh})
algebras.

Let us change the data as follows: suppose that for a nilpotent ideal $N$
of an algebra $A$ the algebra $A/N$ contains a subalgebra
$\overline{S}$ isomorphic to the first Weyl algebra
$\Aw_1=\Bbbk\langle p,q\mid qp-pq=1 \rangle$ or
to the matrix algebra $M_n(\Aw_1)$, $n>1$.
In general, it is not true that the preimage of
$\overline{S}$ in $A$ contains a subalgebra isomorphic to $\Aw_1$ or $M_n(\Aw_1)$,
respectively.
Roughly speaking,
we are going to show that for certain class of algebras an analogue of the
Wedderburn statement holds in these data. The algebras under
consideration are subalgebras of $\End V$ ($V$~is an infinite-dimensional
vector space) satisfying so called TC-condition
(translation-invariance, continuity).
This condition is a new formalization of the class of ordinary algebras associated with
conformal algebras \cite{K1} (c.f. formal distribution algebras \cite{K2}),
so the main results of this paper are
stated for conformal algebras.
In some sense, we continue studying the structure theory of
associative conformal algebras with finite faithful representation
started in \cite{K2,BKL,Ko1}.

The formal definition of a conformal algebra
was introduced in \cite{K1} as an algebraic language describing
the singular part of operator product expansion (OPE) in conformal field
theory. Another important feature of conformal algebras
relates to the notion of a pseudo-tensor category \cite{BD}:
a conformal algebra is just an algebra in the pseudo-tensor category
${\mathcal M}^*(H)$ associated with the polynomial bialgebra $H=\Bbbk[D]$
(see \cite{BDK}). An object of this category is a left unital
$H$-module, and an algebra in ${\mathcal M}^*(H)$
is a module $C\in {\mathcal M}^*(H)$ endowed with $(H\otimes H)$-linear map
$*: C\otimes C \to (H\otimes H)\otimes_H C$.
Associativity, commutativity, and other identities have a natural
interpretation in this language. Note that an ordinary algebra
over a field $\Bbbk $
is just an algebra in ${\mathcal M}^*(\Bbbk)$.

In this context, the obvious way to generalize finite-dimensional
algebras is to consider the class of finite conformal algebras, i.e.,
finitely generated $H$-modules.
Another approach is to consider conformal algebras acting faithfully
on a finitely generated $H$-module, i.e., those
with finite faithful representation.
These conformal algebras are not necessarily finite themselves,
but they are clearly the analogs of subalgebras of $\End U$,
$\dim_{\Bbbk} U<\infty$.
In particular, if $V$ is a finitely
generated $H$-module then the set of all conformal endomorphisms
(see \cite{K2,K1,DK}) forms an associative conformal algebra
denoted by
$\Cend V$, that is a ``conformal analogue'' of $\End U$.
Since $\Cend V$ is an infinite conformal algebra,
the theory of conformal algebras with finite faithful
representation can not be reduced to the theory of finite conformal algebras.
On the other hand, it is unknown whether an arbitrary finite
(Lie or associative) conformal
algebra has a finite faithful representation.

In \cite{K2,DK}, simple and semisimple finite
Lie and associative conformal algebras were described. In the associative
case \cite{K2}, a finite semisimple conformal algebra
is isomorphic to a finite direct sum of current conformal algebras over
simple finite-dimensional algebras. It was shown
in \cite{Z2} that an arbitrary finite associative conformal algebra
$C$ could be presented as $C=S\oplus R$, where $R$ is the maximal nilpotent
ideal of $C$ and $S$ is a semisimple subalgebra isomorphic to $C/R$.

For finite Lie conformal algebras, the similar statement (analogous
to the Levi theorem) is not true: it follows from the description
of cohomologies of the Virasoro conformal algebra \cite{BKV} (see also
\cite{BDK}).

In \cite{Ko1}, the classification of simple and semisimple (more precisely,
semiprime) associative conformal algebras with finite faithful representation
was obtained (see Theorem~\ref{thm-2.1} below). Moreover, any subalgebra $C$
of $\Cend V$ contains a maximal nilpotent ideal $R=\Rad(C)$, and
$C/R$ also has a finite faithful representation. There is a natural problem:
whether $C$ could be presented as $C = S\oplus R$, $S\simeq C/R$?
In this paper we show that the answer is positive if $C/R$ is a unital
conformal algebra, but in the general case the answer is negative.

Throughout this paper, $\Bbbk $ is an algebraically closed field of zero characteristic,
${\mathbb Z}_+$ is the set of non-negative integers. By $H$ we will denote
the polynomial algebra $\Bbbk[D]$ endowed with the $D$-adic topology, i.e.,
the family of basic neighborhoods of zero is given by principal ideals
$(D^n)$, $n\in {\mathbb Z}_+$.

\section{Preliminaries on conformal algebras}

\subsection{Conformal linear maps}

Let $V_1$, $V_2$ be left (unital) $H$-modules, and let
$T_i \in \End V_i$ represent $D$ in $V_i$, $i=1,2$.
A $\Bbbk$-linear function
\[
a: H \to \Hom (V_1,V_2)
\]
is said to be a {\em conformal homomorphism\/}
from $V_1$ to $V_2$ \cite{K2,K1} if:
\begin{itemize}
\item[(i)]
$a$ is continuous
with respect to the $D$-adic topology on $H$ and the finite topology
(see, e.g., \cite{J}) on $\Hom (V_1,V_2)$;
 \item[(ii)]
$a$ is translation-invariant with respect to the operators $T_1$,
$T_2$, i.e., $a(f)T_1 = T_2a(f) + a(f')$ for any $f\in H$
(hereinafter, $f'=\partial_D f$ is the usual derivative).
\end{itemize}
Hereinafter, $\Hom (U,V)$ denotes the set of all $\Bbbk $-linear
maps from $U$ to $V$.

\begin{remark}
It is clear that if $\dim V_1<\infty $,
then the only conformal homomorphism from $V_1$ to $V_2$ is zero.
Therefore, the exponential $e^{\alpha T_1}$, $\alpha\in \Bbbk$, is usually  undefined. However,
one may always interpret (ii) as follows:
\[
e^{-\alpha T_2} a(f(D)) e^{\alpha T_1} = a(f(D+\alpha)), \quad f\in H.
\]
\end{remark}

Let $\Chom(V_1,V_2)$ stands for
the set of all conformal homomorphisms from $V_1$ to $V_2$.
We will write $a(n)$ for
$a(D^n)$ and $a\oo{n} u$ for $a(n)u$,
$a\in \Chom(V_1,V_2)$, $n\in {\mathbb Z}_+$, $u\in V_1$.
By definition, we have
\begin{eqnarray}
& a\oo{n} u = 0 \quad \hbox{for $n$ sufficiently large}, 
 \label{Loc-mod} \\
& a\oo{n}T_1u = T_2(a\oo{n}u) + n a\oo{n-1} u.
  \label{Sesq-mod}
\end{eqnarray}
These properties are equivalent to (i) and (ii) above.

If $V_1 = V_2 = V$ then
$\Chom(V_1,V_2)$ is denoted by $\Cend V$.

One can define the structure of a left $H$-module on $\Chom(V_1,V_2)$
by
\begin{equation}\label{eq-conf-D}
(D a)(f) = -a(f'), \quad f\in H.
\end{equation}

If we have three $H$-modules $V_1$, $V_2$, $V_3$,
then for any $n\in {\mathbb Z}_+$ an operation
\begin{equation}
(\cdot \oo{n} \cdot): \Chom(V_2,V_3)\otimes \Chom( V_1, V_2)
   \to \Chom(V_1,V_3)
\end{equation}
could be defined by the rule
\begin{equation}\label{eq-conf-n}
(a\oo{n} b)(m) = \sum\limits_{s\ge 0} (-1)^s \binom{n}{s} a(n-s)b(m+s),
\quad m\in {\mathbb Z}_+.
\end{equation}
It is easy to see that
\begin{equation}\label{ax-sesq}
Da \oo{n} b = -n a\oo{n-1} b,
\quad
a\oo{n} Db = D(a\oo{n} b) + na\oo{n-1} b.
\end{equation}
Moreover, if $V_1$ is a finitely generated
$H$-module then
for any $a\in \Chom(V_2,V_3)$, $b\in \Chom(V_1,V_2)$
\begin{equation}\label{ax-loc}
a\oo{n} b = 0 \quad \hbox{for $n$ sufficiently large}
\end{equation}
(see \cite{K2,K1}).

\begin{definition}[\cite{K2,K1}]\label{defn-conf}
A {\em conformal algebra\/} $C$ is a left $H$-module
endowed with a family of $\Bbbk$-linear maps
\[
(\cdot \oo{n} \cdot): C\otimes C \to C, \quad n\in {\mathbb Z}_+,
\]
satisfying (\ref{ax-sesq}), (\ref{ax-loc}).
\end{definition}

Axioms (\ref{ax-sesq}) and (\ref{ax-loc}) are known as
{\em sesqui-linearity\/} and {\em locality}, respectively.
The function $N: C\times C \to {\mathbb Z}_+$ given by
\[
N(a,b) = \min\{n\ge 0\mid a\oo{m}b = 0\ \hbox{for all}\ m\ge n\},
\quad  a,b\in C,
\]
is called the {\em locality function\/} on $C$.

If $C$ is a finitely generated $H$-module, then
$C$ is said to be a conformal algebra {\em of finite type}
(or {\em finite\/} conformal algebra).

A conformal algebra $C$ is called associative if \cite{K2,Ro1}
\begin{equation}\label{conf-ass}
(a\oo{n} b) \oo{m} c = \sum\limits_{s\ge 0}
  (-1)^s \binom{n}{s} a\oo{n-s} (b\oo{m+s} c)
\end{equation}
for any $a,b,c\in C$, $n,m\in {\mathbb Z}_+$.

For two subsets $A$, $B$ of a conformal
algebra $C$, denote
\[
A \oo{\omega} B =
 \bigg\{\sum_i a_i\oo{n_i} b_i \mid a_i\in A,\,
  b_i\in B,\, n_i\in {\mathbb Z}_+\bigg\}.
\]
The same notation will be used in other cases when we will work with
a family of binary operations $(\cdot \oo{n}\cdot)$,
$n\in {\mathbb Z}_+$.

It is clear how to define what is a subalgebra or left/right ideal
of a conformal algebra $C$, and what does it mean that $C$
(or an ideal of $C$)
is simple, nilpotent, solvable, prime or semiprime.
In particular, $C$ is semiprime if for any non-zero ideal $I$ of $C$
we have $I\oo{\omega } I \ne 0$. Such conformal algebras
are often called semisimple in the literature
(see, e.g., \cite{K2,Ko1,K1,DK}).
In this paper, we will also use this terminology.

\begin{definition}[\cite{K2,CK}]
Let $C$ be an associative conformal algebra, and let
$V$ be an $H$-module.
A {\em representation} of $C$ on $V$ is an $H$-linear map
\[
\rho : C \to \Cend V
\]
such that
$\rho(a\oo{n}b)= \rho(a)\oo{n}\rho(b)$
for any $a,b\in C$, $n\in {\mathbb Z}_+$.
If $C$ has a representation on $V$ then $V$
is called a (conformal) $C$-module;
if $V$ is a finitely generated $H$-module
then $\rho $ is said to be a representation
{\em of finite type}
(or {\em finite\/} representation).
\end{definition}

If $V$ is a finitely generated $H$-module then
$\Cend V$ endowed with the operations
(\ref{eq-conf-D}), (\ref{eq-conf-n}) is an associative
conformal algebra.
If an associative conformal algebra $C$ has a
finite faithful (i.e., injective) representation on $V$,
then $C$ could be identified with a conformal subalgebra of
$\Cend V$. Such conformal algebras were studied, for example,
in \cite{K2,BKL,Ko1,Re1}, and they are the main objects
of the present work.

\subsection{Correspondence between subalgebras of $\End V$ and $\Cend V$}

Let us first recall
 the following notation \cite{K2,BFK}.
If $C$ is a conformal algebra, $a,b\in C$, $n\in {\mathbb Z}_+ $,
then
\begin{equation}
\{a\oo{n} b\} = \sum\limits_{s\ge 0} \frac{(-1)^{n+s}}{s!} D^{s}(a\oo{n+s} b).
\end{equation}
Relations (\ref{ax-sesq}) imply
\[
\{a\oo{n} Db\} = -n \{a\oo{n-1} b\}, \quad
\{Da\oo{n} b\} = D\{a\oo{n} b\} + n\{a\oo{n-1} b\}.
\]
Moreover,
if $C$ is associative then
\cite{K2,BFK}:
\begin{gather}
a \oo{n} \{ b \oo{m} c\}  =
 \{(a \oo{n} b) \oo{m} c\};   \label{eq2.2.1} \\
\{a \oo{n} (b \oo{m} c)\}   =
   \sum\limits_{s\ge 0} (-1)^s\binom{m}{s}
   \{\{a \oo{m-s} b\} \oo{n+s} c \}  ;
  \label{eq2.2.2}       \\
\{ a \oo{n}\{b \oo{m} c \}\}   =
  \sum\limits_{s\ge 0}(-1)^s \binom{m}{s}
   \{\{ a\oo{n+s} b\} \oo{m-s} c \};
  \label{eq2.2.3} \\
\{a \oo{n} b\} \oo{m} c  =
  \sum\limits_{s\ge 0} (-1)^s \binom{n}{s} a \oo{m+s}(b \oo{n-s} c).
  \label{eq2.2.4}
\end{gather}

Now, let $V$ be a finitely generated $H$-module, $T\in \End V$ represents
the action of~$D$. Consider a correspondence
between conformal subalgebras of $\Cend V$ and ordinary
subalgebras of $\End V$.

For a subset $S\subseteq \End V$
define ${\mathcal F}(S)$ as the set of all $a\in \Cend V$
such that $a(f)\in S$ for all $f\in H$.
For any subset $C \subseteq \Cend V$, define
${\mathcal A}(C)=\{a(f)\mid a\in C,\, f\in H\} \subseteq \End V$.
If $S$ is a subalgebra of $\End V$, then  ${\mathcal F}(S)$
is a conformal subalgebra of $\Cend V$.
If $C$ is a conformal subalgebra of $\Cend V$ then ${\mathcal A}(C)$
is an ordinary subalgebra of $\End V$.
Indeed, it follows from
(\ref{eq-conf-D}) that ${\mathcal A}(C)$ is a subspace, and
(\ref{eq-conf-n}) implies
\begin{equation}\label{oper-prod}
a(n)b(m) = \sum\limits_{s\ge 0} \binom{n}{s} (a\oo{n-s} b)(m+s),
\quad
a,b\in C,\ n,m\in {\mathbb Z}_+.
\end{equation}
Note that
${\mathcal A}({\mathcal F}(S)) \subseteq S$,
$C\subseteq {\mathcal F}({\mathcal A}(C))$.
Denote
$\widetilde{C} = {\mathcal F}({\mathcal A}(C))$.

\begin{lemma}[{\cite[Proposition~3.7]{Ko1}}]\label{lem-2}
Let $C$ be a conformal subalgebra of $\Cend V$,
and let $\{a_n\}_{n\in {\mathbb Z}_+}$ be a sequence of operators
in ${\mathcal A}(C)$ such that $a_n T - T a_n = n a_{n-1}$
for any $n\in {\mathbb Z}_+$. If $\lim\limits_{n\to \infty} a_n = 0$
in the sense of finite topology on $\End V$ then there exists
$a\in \widetilde{C}$ such that $a(n)=a_n$.
\end{lemma}

\begin{lemma}[c.f. {\cite[Proposition~3.10]{Ko1}}]\label{prop-2}
Let $C$ be a conformal subalgebra of $\Cend V$.
\begin{itemize}
\item[\emph{(i)}]
$S={\mathcal A}(C)$ acts on $C$ as on $S$--$S$-bimodule by the rule
\[
a(n)\cdot b = a\oo{n} b,
\quad
b\cdot a(n) = \{b\oo{n} a\}
\]
for
$a,b\in C$, $n\in {\mathbb Z}$;
\item[\emph{(ii)}]
$C$ is a two-sided ideal of $\widetilde{C}$, and
$\widetilde{C}\oo{\omega} \widetilde{C} \subseteq C$.
\end{itemize}
\end{lemma}

\begin{proof}
(i)
It is sufficient to check that the action is correctly
defined. Associativity of the action follows from
(\ref{eq-conf-n}), (\ref{eq2.2.1}), (\ref{eq2.2.3}),
(\ref{oper-prod}).

If $a(n)=0$ for some $a\in \Cend V$,
$n\in {\mathbb Z}_+$, then $a=D^{n+1}x$ for an appropriate
$x\in \Cend V$. Suppose $a(n)=a'(n')$  for some
$a,a'\in C$, $n,n' \in {\mathbb Z}_+$. We may also assume that
$n'\le n$, so that for
\[
x =a - \frac{n'!}{n!}(-D)^{n-n'}a' \in C
\]
we have $x(n) = a(n) - a'(n') =0$, and
$x = D^{n+1}y$.
Therefore,
for any $b\in C$ we have
$0=x \oo{n} b = a\oo{n} b - a'\oo{n'} b$

(ii)
It follows from (i) that
$ \widetilde{C}\oo{\omega }C ,
\{C\oo{\omega} \widetilde{C}\} \subseteq C $.
Locality axiom (\ref{ax-loc}) implies that
$ C \oo{\omega } \widetilde{C} \subseteq C $
as well. Since for any $x\in \widetilde{C}$
 and for any $n\in {\mathbb Z}_+$
there exists $a\in C$ such that $x-a \in D^{n+1}\Cend V$,
we conclude $\widetilde{C}\oo{\omega} \widetilde{C} \subseteq C$.
\end{proof}

\begin{definition}\label{defn-loc-cond}
Let $V$ be a finitely generated $H$-module, and let $S$ be a subalgebra
of $\End V$.
If
$S = {\mathcal A}({\mathcal F}(S))$,
then
we say that $S$ satisfies {\em TC-condition}, or that
$S$ is a {\em TC-subalgebra\/} of $\End V$.
\end{definition}

The most important example appears in the case when $V$ is a
free $H$-module.
Consider the space $V_n = \Bbbk[t]\otimes \Bbbk^n$
as an $H$-module with respect to the representation $D\mapsto T$,
$T: f\otimes u \mapsto tf\otimes u$, $f \in \Bbbk[t]$, $u\in \Bbbk ^n$.
The maximal TC-subalgebra of $\End V_n$
is given by
\[
{\mathcal A}({\mathcal F}(\End V_n))
 = W\otimes \End \Bbbk^n \simeq M_n(W),
\]
where $W$ the subalgebra of $\End \Bbbk[t]$ generated by
$p: f\mapsto tf$ and $q: f\mapsto f'$, $f\in \Bbbk[t]$
(i.e., $W$ is isomorphic to
the first Weyl algebra~$\Aw_1$).

Denote $\Cend V_n$ by $\Cend_n$. For a fixed basis of $V_n$ over $H$,
there is an isomorphism between conformal algebras $\Cend_n$ and
$M_n(\Bbbk[D,v])\simeq  H\otimes M_n(\Bbbk[v])$
\cite{K2,Ko1,Re1}.
Here $v$ is just a formal variable,
an element
\[
a = \sum\limits_{s= 0}^m
\frac{(-D)^s}{s!}\otimes A_s(v) \in H\otimes M_n(\Bbbk[v])
\]
acts on $V_n$ as
\[
a(k) = \sum\limits_{s= 0}^m \binom{k}{s} A_s(p)q^{k-s} \in M_n(W),\quad k\in {\mathbb Z}_+,
\]
and the operations (\ref{eq-conf-n}) are given by
\[
(1\otimes A)\oo{k}(1\otimes B) = 1\otimes A \partial_v^k(B),
\]
$A,B\in M_n(\Bbbk[v])$. From now on, we will identify $\Cend_n$
with $M_n(\Bbbk[D,v])$.

\begin{example}\label{exmp-cur}
The set of matrices $S_0 = M_n(\Bbbk[q])\subset \End V_n$
is a TC-subalgebra. The corresponding conformal algebra
${\mathcal F}(S_0)=M_n(\Bbbk[D])$ is denoted by
$\Curr_n$.
\end{example}

\begin{example}\label{exmp-cend-Q}
If $Q(p)\in M_n(\Bbbk[p])$, then
$W_{n,Q} = M_n(W)Q(p)\subseteq \End V_n$
is a TC-subalgebra, and
$\Cend_{n,Q} = {\mathcal F}(W_{n,Q})=\Cend_n Q(v-D)$
is a conformal subalgebra (even a left ideal) of $\Cend_n$.
If $\det Q(p) \ne 0$ then $\Cend_{n,Q}$ is a simple conformal algebra \cite{BKL}.
\end{example}

\begin{theorem}[\cite{Ko1}]\label{thm-2.1}
Let $C$ be an associative conformal algebra with a finite
faithful representation.
If $C$ is simple, then $C$ is isomorphic either to
$\Curr_n$ or to
$\Cend_{n,Q}$, $n\ge 1$, $\det Q\ne 0$.
If $C$ is semisimple,
then $C$ is a finite direct sum of simple ones.
\end{theorem}

It was also shown in \cite{Ko1} that an arbitrary
associative conformal algebra $C$ with a finite faithful representation
has the maximal nilpotent ideal (radical) $\Rad(C)$,
although $C$ is not necessarily Noetherian.
The following proposition gathers
additional information on the structure
of such conformal algebras.

\begin{proposition}\label{prop-2.2}
Let $C$ be an associative conformal algebra with a finite
faithful representation, and let $R=\Rad (C)$. Then
\begin{itemize}
\item[\emph{(i)}]
 $C/R$ has a finite faithful representation;
\item[\emph{(ii)}]
there exist a finite number of prime ideals of $C$.
The intersection of these ideals is equal to $R$.
\end{itemize}
\end{proposition}

\begin{proof}
(i) Suppose $C\subseteq \Cend V$, $\rank V<\infty$,
$N={\mathcal A}(R)\subseteq \End V$.
There exists $m\ge 1$ such that $N^m=0$.
Then the finitely generated $H$-module
\[
U = V/NV \oplus NV/N^2 V \oplus
  \dots \oplus N^{m-2} V/N^{m-1} V \oplus N^{m-1} V
\]
is also a $C$-module. Since $R$ is a maximal nilpotent ideal,
the annihilator of $U$ in $C$ coincides with $R$, so $U$ is a
faithful $C/R$-module.

(ii) This is a general fact that the prime radical of
an associative conformal algebra $C$
(defined as the minimal nil-ideal $B$ such that $C/B$ has no non-zero
nilpotent ideals)
  is equal to the intersection
of prime ideals: one may follow the proof from \cite{F}
slightly adjusted for conformal algebras. By Theorem~\ref{thm-2.1},
$B=R$ and  $C/R$ has a finite number of ideals,
so there exist only a finite number of prime ideals of $C$.
\end{proof}

\section{Radical splitting problem}\label{sec-3}

\subsection{Lifting of special elements}\label{sec-3.1}

An element $e$ of a conformal algebra $C$ is called
an idempotent if $e\oo{n} e = \delta_{n,0}e$, $n\in {\mathbb Z}_+$
\cite{Z2}.
Two idempotents $e_1,e_2 \in C$ are mutually orthogonal if
$e_1\oo{\omega} e_2 = e_2\oo{\omega } e_1 = 0$.
An idempotent $e\in C$ is said to be a (conformal) unit if
$e\oo{0} x = x$ for any $x\in C$ \cite{Re1}.
For example, the conformal algebra
$\Curr_n\simeq M_n(\Bbbk[D])\subset \Cend_n$
contains a (canonical) unit corresponding to the identity matrix
$\Idd_n \in M_n(\Bbbk[D])$.
In general, a conformal unit is not unique: e.g.,
any element of the form
$Q^{-1}(v)Q(v-D)$,
$\det Q \in \Bbbk\setminus \{0\}$,
 is a unit of
$\Cend_n$.

The structure of unital associative conformal algebras was
considered in details in \cite{Re1,Re2}. Unfortunately,
it is not clear how to join a unit to an arbitrary conformal algebra.

Throughout this section,
$C$ is an associative conformal algebra, $I$ is a nilpotent ideal of $C$.

\begin{lemma}\label{lem-3.1}
\
\begin{itemize}
\item[\emph{(i)}]
If $\bar e\in C/I$ satisfies the condition
$\bar e \oo{0} \bar e = \bar e$
then there exists
$e\in C$
such that
$e+I =\bar e$,
$e\oo{0}e = e$.
\item[\emph{(ii)}]
Suppose that $C$ contains a unit $e_0$ and there exist 
$\bar e_1,\dots,\bar e_N \in C/I$ satisfying the conditions
$\bar e_i \oo{0} \bar e_i = \bar e_i$,
$\{\bar e_i\oo{0}\bar e_0\} = \bar e_i$,
$\bar e_i\oo{n} \bar e_j = 0$ for $i\ne j$, $n\ge 0$,
$i,j=1,\dots,N$.
Then there exist $e_1,\dots, e_N \in C$ such that
$e_i+I = \bar e_i$, $\{e_i\oo{0} e_0\} = e_i$, $e_i\oo{0} e_i = e_i$,
$e_i\oo{0} e_j = 0$ for $i\ne j$,
$i,j=1,\dots , N$.
\end{itemize}
\end{lemma}

\begin{proof}
(i) This statement was proved in \cite{Z2}. The idea is similar to the
lifting of idempotents in ordinary algebras, see, e.g., \cite{Fe}.

(ii) Let $N=1$.
If $e_1'$ is a preimage of $\bar e_1$ given by (i),
then $e_1 = \{e_1'\oo{0} e_0\}$ satisfies the required conditions.

Assume the statement holds for some $N\ge 1$, and let
we are given $N+1$ elements $\bar e_1, \dots, \bar e_N, \bar e_{N+1}\in C/I$
as above.
Suppose we have found $e_1,\dots,e_N\in C$, satisfying the required conditions.
Consider the element
$f = e_0 - e_1 - \dots -e_N\in C$
and the subalgebra
$C_0 = f\oo{0} \{C \oo{0} f\} \subseteq C$
with the nilpotent ideal
$I_0 = C_0\cap I$.
Note that
$\{f\oo{0} e_0\}=f$,
$f \oo{0} e_i = 0$,
$i=1,\dots,N$.
Since $\bar e_{N+1} \in C_0/I_0 \subseteq C/I$,
there exists
$e_{N+1} \in C_0$
such that
$e_{N+1}\oo{0} e_{N+1} = e_{N+1}$.
Presentation
$e_{N+1} = f\oo{0} \{ y\oo{0} f\}$, $y\in C$,
shows that
$\{e_{N+1}\oo{0} e_0\} = e_{N+1}$
and
$e_{N+1}\oo{0} e_i = 0$,
$i=1,\dots,N$.
Moreover, for any
$a\in C$ we have
$(e_i\oo{0} f)\oo{0} a = (e_i\oo{0} e_0 - e_i)\oo{0}(e_0 \oo{0} a) = 0$.
Therefore,
$e_i\oo{0} e_{N+1} = 0$ for $i=1,\dots, N$.
\end{proof}

\begin{lemma}[c.f. \cite{Z2}]\label{lem-3.2}
\

\begin{itemize}
\item[\emph{(i)}]
 If $\bar e\in C/I$ is an idempotent then there exists
an idempotent $e\in C$ such that $e+I=\bar e$.
\item[\emph{(ii)}]
Suppose that $C$ contains a unit $e_0$
and there exists a family
of pairwise mutually orthogonal idempotents
$\bar e_1, \dots, \bar e_N\in C/I$ 
(i.e., $\bar e_i \oo{n} \bar e_j = \delta_{n,0}\delta_{i,j}\bar e_j$)
such that
$\{\bar e_i\oo{0} \bar e_0\} = \bar e_i$.
Then there exist pairwise mutually orthogonal idempotents
$e_1,\dots , e_N \in C$ such that $e_i + I = \bar e_i$.
\end{itemize}
\end{lemma}

\begin{proof}
(i) It was shown in \cite{Z2}.
One may also use the idea of the proof of Proposition~\ref{prop-3.1}
below.

(ii) First, let us find $f_i\in C$, $i=1,\dots, N$, by Lemma \ref{lem-3.1}(ii):
$f_i\oo{0} f_j = \delta_{i,j} f_j$,
$\{f_i\oo{0} e_0\} =f_i$,
$f_i+I = \bar e_i$.
Note that for any $i=1,\dots,N$  subalgebra
$C_i = f_i\oo{0} \{C\oo{0} f_i\}$ contains a preimage of $\bar e_i$.
For every $i=1,\dots, N$, apply (i) to the algebra
$C_i$ with the nilpotent ideal $I\cap C_i$
and find idempotents $e_i\in C_i$, $e_i\in f_i+ I$.
Present $e_k = f_k\oo{0}\{e'_k \oo{0} f_k\}$,
$k=1,\dots, N$.
It is easy to see that
$f_i\oo{0} e_j = f_i\oo{0} f_j \oo{0} \{e'_j \oo{0} f_j\} =0$
and
$e_i\oo{n} e_j = f_i\oo{0}(e_i'\oo{n} (f_i\oo{0} e_j))=0$
for $i\ne j$.
\end{proof}

\begin{lemma}\label{lem-3.3}
Let $C$ be a conformal algebra with a unit $e$,
and assume that there exists $\bar x \in C/I$
such that
$\bar x \oo{0}\bar e = \bar x$,
$\bar e\oo{1}\bar x = \bar e$.
Then there exists a preimage $x\in C$ of $\bar x$ such that
$x\oo{0}e = x$, $e\oo{1} x = e$.
\end{lemma}

\begin{proof}
Let $x_0 \in C$ be a preimage of $\bar x$.
Without loss of generality we may assume that
$x_0 \oo{0} e = x_0$ and $I^2 = 0$. Suppose that
$N(e,x_0)\ge 3$,
where $N(\cdot,\cdot)$ is the locality function on $C$.
Consider
\[
x_{1} = x_0 - \frac{1}{n}(x_0\oo{1}x_0 - x_0),
\quad n=N(e,x_0)-1\ge 2.
\]
It is straightforward to check that $x_1\oo{0}e = x_1$
and $N(e,x_1)\le n$. Indeed,
\begin{multline}\label{eq-3.1.1}
e\oo{n} x_1 = \frac{n+1}{n} e \oo{n}x_0 - \frac{1}{n}e\oo{n}(x_0 \oo{1} x_0) \\
= \frac{n+1}{n} e \oo{n}x_0 - \frac{1}{n}(e\oo{n}x_0)\oo{1}x_0   \\
-\frac{1}{n}\sum\limits_{s=1}^{n-2} \binom{n}{s} (e\oo{n-s}x_0)\oo{s+1} x_0
-(e\oo{1} x_0)\oo{n}x_0.
\end{multline}
Since $(e\oo{m}x_0)\oo{k}x_0 = (e\oo{m}x_0)\oo{0}(e\oo{k} x_0)$
and $e\oo{1}x_0 = e + a$, $a\in I$,
we conclude that the right-hand side of (\ref{eq-3.1.1})
is equal to zero. In the same way, one may show that
$e\oo{m}x_1 = 0$ for any $m>n$.

Therefore, if we choose a preimage $y\in C$ of $\bar x$ such that
$y\oo{0} e = y$ and $N(e,y)$ is minimal, then $N(e,y)=2$.
Suppose
$e\oo{1}y = e+b$, $b\in I$. Since
$e\oo{2}y = e\oo{1}(e\oo{1}y) = e\oo{1}e + e\oo{1}b = 0$,
we have $e\oo{1}b = y\oo{1}b=0$. Moreover, $b\oo{0}e=b$.
Hence,
\[
(e\oo{1}y)\oo{0} (e-b) = e\oo{1}(y-y\oo{0}b), \quad
(e+b)\oo{0}(e-b) = e,
\]
and $x=y-y\oo{0}b$
satisfies the conditions $e\oo{1}x=e$, $x\oo{0}e = x$.

If $I^\nu =0$ for $\nu>2$,
then the lifting of $x$ could be done by induction, using the
sequence
\[
C \to C/I^{\nu-1} \to C/I^{\nu-2} \to \dots \to C/I^2 \to C/I,
\]
as usual.
\end{proof}

\subsection{Unital case: splitting of radical}\label{sec-3.2}

Let us fix a finitely generated $H$-module $V$,
a conformal subalgebra $C$ of $\Cend V$, a nilpotent ideal $I$ of~$C$,
and let $R=\Rad(C)$ be the maximal nilpotent ideal of~$C$.

\begin{proposition}\label{prop-3.1}
Assume $C/I$ contains a subalgebra $\overline{S}$ isomorphic
to $\Curr_N$ for some $N\ge 1$. Then
the preimage of $\overline{S}$ in $C$ contains
a subalgebra $S$ isomorphic to $\Curr_N$.
\end{proposition}

\begin{proof}
Let $\bar e$ be the canonical unit of $\overline{S}$. By Lemma~\ref{lem-3.2}(i)
there exists an idempotent $e\in C$ which is a preimage of $\bar e$.
If $N=1$ then the $H$-span of $e$ is isomorphic to $\Curr_1$.
If $N>1$ then consider the subalgebra $C_e = e\oo{0}\{C\oo{0}e\} \subseteq C$
which is unital. Moreover, $\overline{S}\subseteq C_e/I\cap C_e$,
so $C_e/I\cap C_e$ contains a family of
pairwise mutually
orthogonal idempotents
$\bar e_1, \dots, \bar e_N \in \overline{S}$
corresponding to diagonal matrix units of $M_N(\Bbbk)$.
We may apply Lemma~\ref{lem-3.2}(ii)
to find orthogonal idempotents $e_1,\dots, e_N \in C_e$.

Our aim is to build a system of (conformal) matrix units in $C$,
i.e., a family of $e_{ij}\in C$, $i,j=1,\dots,N$,
such that
\begin{equation}
e_{ij}\oo{n} e_{kl} = \delta_{n,0}\delta_{j,k} e_{il}, \quad n\ge 0.
\end{equation}
Let us choose some preimages $v_{1j}, v_{i1}\in C_e$, $i,j=2,\dots,N$,
of the corresponding matrix units.
We may assume $e_1\oo{0}v_{1j}\oo{0} e_j = v_{1j}$,
$e_i \oo{0} v_{i1} \oo{0} e_1 = v_{i1}$.
Since
\[
v_{1i}\oo{0} v_{i1} = e_1 + a_i, \quad a_i\in I\cap C_e,
\]
and $a_i$ is nilpotent with respect to the 0-product, we may find
$b_i= -a_i + a_i\oo{0}a_i - a_i\oo{0}a_i\oo{0}a_i + \dots \in I\cap C_e$
such that $a_i + b_i + a_i\oo{0} b_i = 0$ and $e_1\oo{0}b_i\oo{0} e_1 = b_i$.
Then the elements
\begin{equation} \label{eq-3.4.1}
f_{1j} = v_{1j}, \quad f_{i1} = v_{i1}+v_{i1}\oo{0} b_i, \quad i,j=2,\dots,N,
\end{equation}
satisfy the following relations:
\begin{eqnarray}
&  e_1\oo{n} e_1 = \delta_{n,0}e_1,
  \quad
     e_1 \oo{0} f_{1j} = f_{1j}, \quad f_{i1}\oo{0} e_1 = f_{i1}, \\
&  f_{1j}\oo{n} f_{i1} = 0 \ \hbox{for}\ i\ne j,\  n\ge 0,
 \quad  f_{1i}\oo{0} f_{i1} = e_1, \\
&  f_{1j}\oo{n} f_{1i} = f_{i1}\oo{n} f_{j1} = 0, \  n\ge 0\\
&  f_{i1}\oo{n} e_1 = f_{1j}\oo{n} e_j = 0 \ \hbox{for}\ n\ge 1.
\end{eqnarray}

Denote $V_j = f_j\oo{0} V \subseteq V$,
$C_{ij} = e_i\oo{0}\{C \oo{0} e_j\}$,
$i,j = 1,\dots, N$. It is easy to see that
$C_{ij}\subseteq \Chom(V_j,V_i)$, $e_i\in C_{ii}$,
$f_{1j}\in C_{1j}$, $f_{i1}\in C_{i1}$.
Let
$S_{ij} = {\mathcal A}(C_{ij})\equiv \{a(n) \mid a\in C_{ij}, \, n\ge 0\}
 \subseteq \Hom (V_j,V_i)$.
Note that $S_{ii}$ is a TC-subalgebra of $\End V_i$
and $S_{ij}$ is an $S_{ii}$--$S_{jj}$-bimodule. Moreover,
$S_{ij}S_{jk} \subseteq S_{ik}$.

Let us fix $j\in \{2,\dots, N\}$ and consider
the sequence of operators
\begin{equation}
b(n) = (f_{j1}(1)f_{1j}(0))^n \in S_{jj}, \quad n\ge 0.
\end{equation}
For any  $n\ge 0$
there exists $d_n\in I\cap C_{jj}$
such that $b(n)-e_j(n) = d_n(n)$.
Therefore,
\[
(b(n) - e_j(n))( (I\cap C_{jj})^k\oo{\omega} V_j)
 \subset  (C_{jj}\cap I)^{k+1}\oo{\omega }V_j.
\]
Since $\lim\limits_{n\to \infty}e_j(n)= 0$ in the sense of finite topology,
for any $u\in V_j$ there exists $n_1\ge 0$ such that
$u_1 = b(n)u \in (C_{jj}\cap I)\oo{\omega } V_j$ for $n\ge n_1$.
In the same way, there exists $n_2\ge 0$ such that
$b(n) u_1 \in (C_{jj}\cap I)^2\oo{\omega} V_j$ for $n\ge n_2$,
and so on. Since $b(n)b(m) = b_{n+m}$, we conclude that
$\lim\limits_{n\to \infty}b(n) = 0$.
Moreover,
$[b(n), D] = nb(n-1)$, so by Lemma~\ref{lem-2}
there exists an element $h_j\in \widetilde{C}_{jj}$
such that $h_j(n)=b(n)$. Since $C_{jj}$ is unital, $h_j \in C_{jj}$
by Proposition~\ref{prop-2}(ii).
Moreover, $(h_j - e_j)(n) = d_n(n)$, so there exists
 $\Delta_n \in C_{jj}$
such that $h_j - e_j = d_n + D^{n+1}\Delta_n$. Hence,
$\{e_j \oo{n} (h_j - e_j)\}\in I$ for any $n\ge 0$, and
$e_j \oo{n}(h_j - e_j) \in I$ because of locality.
In particular, $e_j \oo{0} (e_j - h_j) = e_j -h_j \in I$.

Let $e_{1j} = f_{1j}\oo{0} h_j$. It is clear that $e_{1j} - f_{1j}\in I$.
Remind that $e_1 = f_{1j}\oo{0} f_{j1}$ and $e_1(m)=e_1(1)^m$.
Then
\begin{multline}
S_{1j}\ni (e_1 \oo{n} e_{1j})(m)
=
\sum\limits_{s\ge 0} (-1)^s \binom{n}{s} e_1(n-s) e_{1j}(m+s) \\
=
\sum\limits_{s\ge 0} (-1)^s \binom{n}{s}
   (f_{1j}(0)f_{j1}(1))^{n-s}f_{1j}(0) (f_{j1}(1)f_{1j}(0))^{m+s} \\
=
\sum\limits_{s\ge 0} (-1)^s \binom{n}{s}
f_{1j}(0) (f_{j1}(1)f_{1j}(0))^{n+m}
=
\begin{cases}
0, & n\ge 1;\\
e_{1j}(m), & n=0.
\end{cases}
\end{multline}
Hence, $e_1\oo{n} e_{1j} = \delta_{n,0}e_{1j}$.
Moreover, since $f_{j1}(m) = f_{j1}(0)e_1(m)$, we have
\begin{multline}
(e_{1j}\oo{n} f_{j1})(m)
=
 \sum\limits_{s\ge 0} \binom{n}{s} e_{1j}(n-s)f_{j1}(m+s) \\
=
 \sum\limits_{s\ge 0} \binom{n}{s} f_{1j}(0)(f_{j1}(1)f_{1j}(0))^{n-s}
 f_{j1}(0) (f_{1j}(0)f_{j1}(1))^{m+s}                \\
=
 \sum\limits_{s\ge 0} \binom{n}{s} (f_{1j}(0)f_{j1}(1))^{n+m}
=
\begin{cases}
0, & n\ge 1; \\
e_1(m), & n=0.
\end{cases}
\end{multline}
Therefore,
$e_{1j}\oo{n} f_{j1} = \delta_{n,0}e_1$.

Now it is clear that
\[
e_{1j},\quad e_{ij}= f_{i1}\oo{0} e_{1j}, \quad e_{i1} = f_{i1},
\quad i,j=2,\dots, N,
\]
is a system of matrix units in $C$.
The $H$-linear span of $\{e_{ij}\mid i,j=1,\dots,N\}$
is a conformal algebra $S\subseteq C$ isomorphic to $\overline{S}$.
\end{proof}

\begin{proposition}\label{prop-3.2}
Assume $C/I$ contains a subalgebra $\overline{S}$ isomorphic
to $\Cend_N$ for some $N\ge 1$. Then the preimage of $\overline{S}$
in $C$ contains a subalgebra $S$ isomorphic to $\Cend_N$.
\end{proposition}

\begin{proof}
By Proposition~\ref{prop-3.1}
there exists a subalgebra $S_0\subset C$ isomorphic to $\Curr_N
\subset \Cend_N$.
Let $e_{11}\in S_0$ be the preimage of the corresponding matrix unit
${\mathrm I}_{11}\in \Cend_N$.
The subalgebra $C_1 = e_{11}\oo{0}\{C\oo{0}e_{11}\}$ is unital and contains
a preimage of the
element $\bar x_1\in C/I$ corresponding to $v{\mathrm I}_{11} \in \Cend_N$.
By Lemma~\ref{lem-3.3} there exists $x_1 \in C_1$ such that
$e_{11}\oo{1} x_1 = e_{11}$, $x_1\oo{0} e_{11} = e_{11}$.
Denote
\[
x = \sum\limits_{i=1}^N e_{i1}\oo{0} x_1 \oo{0} e_{1i} \in C,
\]
where $e_{1i}, e_{i1}$ are the preimages of the corresponding matrix
units constructed by Proposition~\ref{prop-3.1}. The element~$x$
is a preimage of
$v\Idd_N \in M_N(\Bbbk[D,v])\simeq \Cend_N$.
It is straightforward to check that
\begin{eqnarray}
&  e_{ij}\oo{0} x = x \oo{0} e_{ij}, \\
&  e_{ij}\oo{1} x = e_{ij}, \quad e_{ij}\oo{n} x = 0, \ n\ge 2, \\
&  x\oo{n}e_{ij} = 0, \ n\ge 1, \\
&  \sum\limits_i e_{ii}\oo{0}x = \sum_i x\oo{0} e_{ii} = x.
\end{eqnarray}
These are the defining relations of the algebra $\Cend_N$.
Since $\Cend_N$ is simple, the conformal algebra generated in $C$ by $S_0\cup\{x\}$
is isomorphic to
$\Cend_N$.
\end{proof}

\begin{theorem}\label{thm-3.1}
Let $C$ be a conformal subalgebra of $\Cend V$, and let $R=\Rad(C)$
be the maximal nilpotent ideal of $C$. If $C/R$ contains a
unit then there exists a semisimple subalgebra $S$ of $C$
such that $C = S\oplus R$.
\end{theorem}

\begin{proof}
It follows from Proposition~\ref{prop-2.2}
that
\[
C/R \simeq \bigoplus\limits_i \overline{C}_i, \quad
\overline{C}_i \simeq \Cend_{n_i} \  \hbox{or}\ \overline{C}_i\simeq \Curr_{n_i}.
\]
Consider the canonical units $\bar e_i$ of $\overline{C}_i$,
then their sum $\bar e = \sum_i \bar e_i$ is a unit of $C/R$.
Let $e\in C$ be the idempotent (preimage of $\bar e$)
constructed by Lemma~\ref{lem-3.2}(i).
The analogue of Pierce decomposition
\begin{multline}\label{eq-3P}
C = e\oo{0} \{C\oo{0} e\} \oplus (1-e)\oo{0} \{C\oo{0} e\} \\
\oplus e\oo{0}\{C\oo{0} (1-e)\}
\oplus (1-e)\oo{0} \{C \oo{0} (1-e)\}
\end{multline}
could be constructed (see \cite{Z2}). Here
$(1-e)\oo{0}a = a-e\oo{0}a$, $\{a\oo{0}(1-e)\} = a - \{a\oo{0}e\}$.
The first summand $C_0 = e\oo{0}\{C\oo{0}e\}$ in (\ref{eq-3P})
is a unital conformal algebra, and
$C_0/R\cap C_0 \simeq C/R$. It is sufficient to show that $C_0$ could be presented
as $S\oplus (R\cap C_0)$.

By Lemma \ref{lem-3.2}(ii), there exist orthogonal idempotents $e_i\in C_0$
such that $e_i + I = \bar e_i$.
If $\overline{C}_i \simeq \Curr_{n_i}$ then the subalgebra
$C_i = e_i\oo{0} \{C_0 \oo{0} e_i\}$ satisfies
the conditions of Proposition~\ref{prop-3.1}.
If $\overline{C}_i\simeq \Cend_{n_i}$ then one may apply Proposition~\ref{prop-3.2}.
In any case, there exist $S_i \subseteq C_i \subseteq C_0$ such that
$S_i\simeq \overline{C}_i$. The direct sum $S=\bigoplus_i S_i$ is a subalgebra
of $C_0$ isomorphic to $C/R$. Therefore,
$C_0 = S\oplus (C_0\cap R)$ and $C = S\oplus R$.
\end{proof}

\begin{remark}
Although it is unknown whether an arbitrary finite associative conformal algebra
is a subalgebra of $\Cend V$ for a finitely generated $H$-module $V$,
it is possible to derive the radical splitting theorem for finite conformal algebras
\cite{Z2} from Theorem~\ref{thm-3.1}. Indeed, if $C$ is a finite associative
algebra, $I=\{a\in C \mid a\oo{\omega } C = 0\}$, then $C/I \subseteq \Cend C$.
Since $I\oo{\omega }C =0$, it is easy to raise the semisimple part of $C/I$
into~$C$.
\end{remark}

\begin{corollary}\label{cor-5}
Let $V$ be a finitely generated $H$-module, and let
$A\subseteq \End V$ be a TC-subalgebra.
\begin{itemize}
 \item[\emph{(i)}]
If $A$ is semiprime then
$A$ is a finite direct sum
of algebras
isomorphic either to
$M_{n}(\Bbbk[q])$ or $M_{n}(\Aw_1)Q(p)$, $Q=\diag(f_1,\dots,f_{n})$
is a matrix in the canonical diagonal form, $\det Q \ne 0$.
 \item[\emph{(ii)}]
The Jacobson radical $J(A)$ of $A$ is nilpotent. The
semisimple image $A/J(A)$ is also a TC-subalgebra of $\End U$ for
an appropriate finitely generated $H$-module $U$.
 \item[\emph{(iii)}]
If $A/J(A)$ contains a unit then there exists a
semisimple TC-subalgebra $B\subseteq \End V$
such that $A = B\oplus J(A)$.
\end{itemize}
\end{corollary}

\begin{proof}
Statement (i) follows directly from Theorem~\ref{thm-2.1}
and Definition~\ref{defn-loc-cond}:
if $A$ is semiprime then ${\mathcal F}(A)\subseteq \Cend V$ is semiprime.

(ii) Let $C={\mathcal F}(A)$, $R=\Rad(C)$, $N={\mathcal A}(R)\subseteq A$.
Since $N$ is a nilpotent ideal, $N\subseteq J(A)$.
Let $U$ be the faithful finite $C/R$-module built by
Proposition~\ref{prop-2.2}(i). Then $C/R \subseteq \Cend U$
and ${\mathcal A}(C/R)\subseteq \End U$ is a TC-subalgebra.
It is straightforward to show that
\[
A/N \to {\mathcal A}(C/R), \quad a(n)+N \mapsto (a+R)(n), \quad a\in C,\ n\in {\mathbb Z}_+,
\]
is a well-defined isomorphism of algebras. Thus, $A/N$ is semisimple and $J(A) = N$.

Now (iii) easily follows from Theorem~\ref{thm-3.1}.
\end{proof}

\subsection{General case: a counterexample}\label{sec-3.3}

Let us show that the unitality condition
in Theorem~\ref{thm-3.1} is essential.
Consider the following subset of $\Cend_2$:
\begin{equation}
C = \Bbbk[v-D] \left\{ a(f,g) \mid
   f(v), g(v)\in \Bbbk[v] \right\}.
\end{equation}
where
\[
a(f,g) = \begin{pmatrix}
                     v^2f(v)\quad & v^2f(v) + v^2g(v)(v-D)^2 \\
                     0\quad       & f(v)(v-D)^2
                        \end{pmatrix}
\]
It is straightforward to check that
\begin{multline}                                       \label{eqn-4.1}
a(f_1,g_1) \oo{n} a(f_2, g_2) \\
= a(f_1(v^2f_2)^{(n)},
   f_1(v^2g_2)^{(n)} + f_1f_2^{(n)} + g_1(v^2f_2)^{(n)}).
\end{multline}
Moreover, for any $h\in \Bbbk[v]$ we have
$ h(v)a(f,g) = a(hf,hg)$, thus, $DC\subseteq C$.
Therefore, $C$ is a conformal subalgebra of $\Cend_2$, and $\Bbbk[v]C = C$.

\begin{lemma}\label{lem-4.1}
$\Rad (C) = \Bbbk[v-D]\{a(0,g) \mid g\in \Bbbk[v]\}$.
\end{lemma}

\begin{proof}
Relation ``$\supseteq$'' is obvious. The converse is also true:
assuming
$a = \sum_k (v-D)^k a(f_k, g_k) \in \Rad(C)$ we obtain
that $\sum_k (v-D)^k v^2f_k(v) \in \Cend_1$ is nilpotent.
But $\Cend_1$ contains no non-zero nilpotent elements,
so $f_k = 0$ for any~$k$.
\end{proof}

\begin{lemma}\label{lem-4.2}
$C/\Rad(C)\simeq \Cend_{1,v^2}$.
\end{lemma}

\begin{proof}
Let us define the map
\[
\theta : C\to \Cend_1, \quad
  \sum\limits_k (v-D)^k a(f_k,g_k)
  \mapsto \sum\limits_{k} f_k(v) (v-D)^{k+2}.
\]
This is a well-defined homomorphism of conformal algebras,
$\Image \theta = \Cend_{1,v^2}$,
and it follows from Lemma \ref{lem-4.1} that
$\Kernel \theta = \Rad(C)$.
\end{proof}

Suppose the radical splitting property holds for the conformal algebra $C$, i.e.,
there exists a conformal subalgebra $S\subseteq C$
such that
$C = S\oplus \Rad(C)$, $S\simeq \Cend_{1,v^2}$.
For any $f(v)(v-D)^2\in \Cend_{1,v^2}$, consider its preimage $x(f)$ in~$S$.
This element can be written as
\begin{equation}					\label{eqn-4.2}
x(f)= \begin{pmatrix}
  v^2f(v)\quad  & v^2f(v) + v^2 g(v,D)(v-D)^2 \\
  0\quad   & f(v)(v-D)^2
\end{pmatrix} ,
\end{equation}
where $g(v,D)$ is uniquely defined by~$f$.
Therefore, we may define a linear map
\[
\psi: \Bbbk[v] \to \Bbbk[v,D], \quad \psi (f) = g.
\]
Since $S$ is assumed to be a subalgebra, $\psi $ has to satisfy the following relation:
\begin{equation}			\label{eqn-4.3}
\psi (f_1 \oo{n} v^2f_2) = f_1\oo{n} f_2 + \psi(f_1)\oo{n} v^2f_2 + f_1\oo{n} v^2\psi(f_2).
\end{equation}
In particular, $\psi $ is completely defined by the value $\psi(1)$.

As a generic polynomial, $\psi(1)$ can be uniquely written as
\begin{equation}			\label{eqn-4.4}
\psi(1) = \sum\limits_{k\ge 0} a_k(v)(v-D)^k.
\end{equation}
Since $\psi(1\oo{2} v^2) = 2\psi(1)$, we have
\[
2\bigg ( \sum\limits_{k\ge 0} a_k(v)(v-D)^k \bigg )
= \sum\limits_{k\ge 0} a_k(v)(v-D)^k \oo{2} v^2
 + 1\oo{2} v^2\bigg ( \sum\limits_{k\ge 0} a_k(v)(v-D)^k \bigg )
\]
or
\[
2a_0(v)+\sum\limits_{k\ge 1} 2a_k(v)(v-D)^k
=
\sum\limits_{k\ge 0} a_k(v)(v^{k+2})''
+
(v^2a_0(v))'' + \sum\limits_{k\ge 1} (v^2a_k(v))'' (v-D)^k.
\]
Comparing terms at $(v-D)^k$, $k\ge 0$, we obtain
\begin{gather}
2a_k(v) = (v^2a_k(v))'', \quad k\ge 1,
							         \label{eqn-4.5-1} \\
2a_0(v) = (v^2a_0(v))'' + \sum\limits_{k\ge 0} a_k(v)(v^{k+2})''.
									\label{eqn-4.5-2}
\end{gather}
It follows from (\ref{eqn-4.5-1}) that $a_k\in \Bbbk $ for $k\ge 1$, and then
(\ref{eqn-4.5-2}) implies
$a_0(v) = -\sum\limits_{k\ge 1} a_kv^k $.
Finally,
\[
\psi(1) = \sum\limits_{k\ge 1} (a_k(v-D)^k - a_kv^k).
\]

Now we may use (\ref{eqn-4.3}) to deduce the explicit formula for $\psi (f)$,
$f\in \Bbbk[v]$:
\begin{equation}			\label{eqn-4.6}
\psi(f) = f(v)\psi(1).
\end{equation}
Let us prove (\ref{eqn-4.6}) by induction on $\deg f$. For $\deg f=0 $ we have done. If
(\ref{eqn-4.6}) is true for some $f(v)$, then it follows from (\ref{eqn-4.3}) that
the same formula is true for
$vf(v) = \frac{1}{2}f\oo{1}v^2$.

But (\ref{eqn-4.6}) does not satisfy (\ref{eqn-4.3}): it is sufficient to consider
\[
\psi(1\oo{1}v^3) = 1\oo{1}v + \psi(1)\oo{1}v^3 + 1\oo{1}v^2\psi(v)= 1 + 3\psi(v^2),
\]
which is not equal to $3\psi(v^2)$.

Therefore, such a subalgebra $S$ does not exist, so Theorem~\ref{thm-3.1} does not hold
in general.

\begin{remark}\label{rem-4.4}
The first summand $f_1\oo{n} f_2$ of the right-hand side of (\ref{eqn-4.3})
plays the role of a Hochschild's 2-cocycle. This is not a cocycle in the sense of
\cite{BKV}, where the basics of the Hochschild theory for conformal algebras have been
introduced. Nevertheless, the analogy with the classical theory
\cite{H45} is clear: we have actually
proved that the map $(f_1,f_2)\to f_1\oo{n}f_2$ is not a ``coboundary''.
\end{remark}

\subsection*{Acknowledgements}
I am grateful to Leonid Bokut and Alexander Pozhidaev for useful
comments, to Efim Zel'manov for his contribution in Proposition
\ref{prop-2.2}(i). I would like to acknowledge Igor L'vov by
mention the joint unpublished result on the structure of the prime radical
of an associative conformal algebra
in Proposition~\ref{prop-2.2}(ii).


\begin{thebibliography}{26}

\bibitem{Az}
G.~Azumaya, On maximally central algebras,
{\it  Nagoya Math. J.} {\bf 2} (1951) 119--150.

\bibitem{Fe}
C.~Feldman, The Wedderburn principal theorem in Banach algebras,
{\it Proc. Amer. Math. Soc.} {\bf 2} (1951) 771--777.

\bibitem{W}
J.~A.~Wehlen,
Splitting properties of extensions of the Wedderburn principal theorem,
in: Rings, extensions, and cohomology (Evanston, IL, 1993), 223--240,
Lecture Notes in Pure and Appl. Math. {\bf 159}
(Dekker, New York, 1994).

\bibitem{Ze}
D.~Zelinsky,
Raising idempotents, {\it Duke Math. J.} {\bf 21} (1954) 315--322.

\bibitem{Al}
A.~A.~Albert, The Wedderburn principal theorem for Jordan algebras,
{\it  Ann. of Math. (2)} {\bf 48} (1947) 1--7.

\bibitem{E}
A.~Elduque, On a class of Malcev superalgebras,
{\it J. Algebra} {\bf 173} (1995) no.~2, 237--252.

\bibitem{Sc1}
R.~D.~Schafer,
The Wedderburn principal theorem for alternative algebras,
{\it  Bull. Amer. Math. Soc.} {\bf 55} (1949) 604--614.

\bibitem{Sc2}
R.~D.~Schafer,
On structurable algebras, {\it J. Algebra} {\bf 92} (1985) no.~2, 400--412.

\bibitem{SvO}
D.~Stefan and F.~Van Oystaeyen,
The Wedderburn-Malcev theorem for comodule algebras,
{\it  Comm. Algebra} {\bf 27} (1999) no.~8, 3569--3581.

\bibitem{Rh}
D.~J.~Rodabaugh,
On the Wedderburn principal theorem,
{\it  Trans. Amer. Math. Soc.} {\bf 138} (1969) 343--361.

\bibitem{K2}
V.~G.~Kac,
Formal distribution algebras and conformal algebras,
in: XII-th International Congress in Mathematical Physics (ICMP'97), Brisbane
(Internat. Press, Cambridge, MA, 1999) pp.~80--97.

\bibitem{BKL}
C.~Boyallian, V.~G.~Kac and J.~I.~Liberati,
On the classification of subalgebras of $\Cend_N$ and $\gc_N$,
{\it J.~Algebra} {\bf 260} (2003) no.~1, 32--63.

\bibitem{Ko1}
P.~S.~Kolesnikov,
Associative conformal algebras with finite faithful representation,
{\it Adv. Math.}, to appear; preprint math.QA/0402330.

\bibitem{K1}
V.~G.~Kac,  {\it Vertex algebras for beginners,}
Univ. Lect. Series {\bf 10}
(AMS, Providence, RI, 1997).

\bibitem{BD}
A.~A.~Beilinson and V.~G.~Drinfeld,
{\it Chiral algebras}, Amer. Math. Soc. Colloquium Publications {\bf 51}
(AMS, Providence, RI, 2004).

\bibitem{BDK}
B.~Bakalov, A.~D'Andrea and V.~G.~Kac,
Theory of finite pseudoalgebras,
{\it Adv. Math.} {\bf 162} (2001) no.~1, 1--140.

\bibitem{DK}
A.~D'Andrea and V.~G.~Kac,
Structure theory of finite conformal algebras,
{\it Sel. Math., New Ser.} {\bf 4} (1998) 377--418.

\bibitem{Z2}
E.~I.~Zel'manov,
Idempotents in conformal algebras,
in: Proc. of Third Internat. Alg. Conf. (Y.~Fong et al, eds) (2003), 257--266.

\bibitem{BKV}
B.~Bakalov, V.~G.~Kac and A.~A.~Voronov, Cohomology of conformal algebras,
{\it Comm. Math. Phys.} {\bf 200} (1999) no.~3, 561--598.

\bibitem{J}
N.~Jacobson, {\it Structure of rings},
Amer. Math. Soc. Colloquium Publications {\bf 37}
(AMS, Providence, RI, 1956).

\bibitem{Ro1}
M.~Roitman,
On free conformal and vertex algebras,
{\it J. Algebra} {\bf 217} (1999) no.~2, 496--527.

\bibitem{CK}
S.-J.~Cheng and V.~G.~Kac,
Conformal modules, {\it Asian J. Math.} {\bf 1} (1997), 181--193.

\bibitem{Re1}
A.~Retakh,
Associative conformal algebras of linear growth,
{\it J. Algebra} {\bf 237} (2001) no.~2, 769--788.

\bibitem{BFK}
L.~A.~Bokut, Y.~Fong, W.-F.~Ke and P.~S.~Kolesnikov,
Gr\"obner and Gr\"obner--Shirshov bases in algebra and conformal algebras
(Russian), {\it Fundam. Prikl. Mat.} {\bf 6} (2000) no.~3, 669--706.

\bibitem{F}
C.~Faith, {\it Algebra. II. Ring theory},
Grundlehren der Mathematischen Wissenschaften, No.~191
(Springer-Verlag, Berlin, New York, 1976).

\bibitem{Re2}
A.~Retakh,
On associative conformal algebras of linear growth II,
{\it J. Algebra}, to appear; preprint math.RA/0212168.

\bibitem{H45}
G.~Hochschild, On the cohomology groups of an associative algebra,
{\it Ann. Math.} {\bf 46} (1945) no.~1, 58--67.


\end{thebibliography}
\end{document}